\documentclass[a4paper,12pt]{amsart}
\usepackage[english]{babel}

\usepackage[textwidth=16.5cm, textheight=24.8cm, left=1.5cm, right=1.5cm]{geometry}

 \usepackage{amsmath,amssymb,amsthm,url}
\usepackage{dsfont,textcomp,graphicx,overpic}
\usepackage[usenames,dvipsnames]{color}
\usepackage{hyperref}
\hypersetup{
   colorlinks   = true, 
   urlcolor     = blue, 
   linkcolor    = blue, 
   citecolor   = red 
}
\usepackage{tikz}
\usetikzlibrary{positioning, calc, intersections, through, arrows}
\usetikzlibrary{graphs}
\usetikzlibrary{graphs.standard}
\usetikzlibrary{patterns}

\newtheoremstyle{mydefinition}
  {\medskipamount}
  {\medskipamount}
  {\normalfont}
  {\parindent}
  {\bfseries}
  {.}
  { }
  {}
 
\long\def\comment#1\endcomment{}

\newcommand{\diag}{\operatorname{diag}}

\def\t{\widetilde}
\def\R{{\mathbb R}} \def\Z{{\mathbb Z}}

\newcommand{\vv}[1]{\overrightarrow{\mathstrut#1}}

\newcommand{\simple}[1]{}
\newcommand{\sophi}[1]{#1}

\theoremstyle{plain}
\newtheorem{theorem}{Theorem}[section]

\theoremstyle{mydefinition}
\newtheorem{pr}[theorem]{}

\begin{document}

\title{Cycles in graphs and in hypergraphs: 
results and problems} 

\author{E. Alkin, S. Dzhenzher, O. Nikitenko,  A. Skopenkov, A. Voropaev}
\thanks{\emph{E. Alkin, S. Dzhenzher,  A. Skopenkov:} Moscow Institute of Physics and Technology. 
\newline
\emph{O. Nikitenko:} Altay Technical University. 
\newline
\emph{A. Skopenkov:} Independent University of Moscow, \url{https://users.mccme.ru/skopenko/}. 
\newline
This text (except for problems marked by star) was presented in 2023 at the International Summer Conferences of the Tournament of Towns, 
see \url{https://en.wikipedia.org/wiki/Tournament_of_the_Towns}.  
We are grateful to I. Bogdanov, A. Ryabichev and O. Styrt for useful discussions, and to A. Ryabichev for preparing computer versions of some figures.} 
 
\date{}

\maketitle

\begin{abstract}
This is an expository paper. A 1-cycle in a graph is a set C of edges such that every vertex is contained in an even number of edges from C. E.g., a cycle in the sense of graph theory is a 1-cycle, but not vice versa. It is easy to check that the sum (modulo 2) of 1-cycles is a 1-cycle. In this text we study the following problems: to find

$\bullet$ the number of all 1-cycles in a given graph;

$\bullet$ a small number of 1-cycles in a given graph such that any 1-cycle is the sum of some of them.

We also consider generalizations (of these problems) to graphs with symmetry, and to 2-cycles in 2-dimensional hypergraphs.
\end{abstract}

\tableofcontents

\section*{Preface}
 
\textbf{Informal description of main results.}

A (simplicial) \textbf{$1$-cycle} (modulo $2$) in a graph is a set $C$ of edges such that every vertex is 
contained in an even number of edges from $C$. 
E.g., a cycle in the sense of graph theory is a $1$-cycle, but not vice versa. 
It is easy to check that the sum (modulo~$2$) of $1$-cycles is a $1$-cycle. 

In this text we study the following problems: to find 

$\bullet$ the number of all 1-cycles in a given graph; 

$\bullet$ a small number of 1-cycles in a given graph such that any 1-cycle is the sum of some of them.
 
We also consider generalizations (of these problems) to graphs with symmetry, and to \emph{$2$-cycles} in \emph{$2$-dimensional hypergraphs} (all the italicized  notions are defined later). 
Main problems are \ref{l:gencom1}.\ref{en:gencom1:rel}, \ref{p:knn}.\ref{en:knn:rel}, 
\ref{pr:orinumber}.\ref{en:orinumber:conngraph}, 
\ref{p:tkn}.\ref{en:tkn:4}\ref{en:tkn:6}, and 
\ref{l:gencom}.cf, \ref{pr:rook}.cf. 
The peaks of this text are results on $1$- and $2$-cycles in \emph{the square of a graph} (Problems \ref{stcycles2}.c, \ref{t:kunn}.cd, \ref{t:kunneth}.ab and \ref{l:h2sym}.d). 

The notion of a 1-cycle and its generalizations have many applications in topology, for simpler ones see e.g. \cite[\S\S 4.11, 6, 9]{Sk20}, \cite[\S\S 1, 8, 9]{Sk}. 
However elementary, this text is motivated by frontline of research, see \cite{FH10, MS17, SS23} and the references therein. 
Open problems are  \ref{stcyclesdel}.cdf and \ref{l:h2sym}.f. 

For an exposition of this material in more details and not as a sequence of problems see \cite{DMN+}.  
 
\medskip
\textbf{Learning by doing problems}

In this text we expose a theory as a sequence of problems, see e.g. \cite[Introduction, Learning by doing problems]{Sk21m} and the references therein. 
Most problems are useful theoretical facts. 
So this text could in principle be read even without solving problems or looking to \S\ref{s:answ} `Answers, hints and solutions'. 

Problems are numbered, the words `problem' are omitted. 
If a mathematical statement is formulated as a problem, then the objective is to prove this statement.
Open-ended questions are called 
{\bf riddles}; here one must come up with a clear wording, and a proof. 
If a problem is named `theorem' (`lemma', `corollary', etc.), then this statement is  considered to be more important.
Usually we \emph{formulate} beautiful or important statement \emph{before} giving a sequence of results (lemmas, assertions, etc.) which constitute its \emph{proof}.
We give hints on that after the statements but we do not want to deprive you of the pleasure of finding the right moment when you finally are ready to prove the statement.
In general, if you are stuck on a certain problem, try looking at the next ones; they may turn out to be helpful.

Important definitions are highlighted in \textbf{bold} for easy navigation.

\medskip
\textbf{On presentation at the Summer conference}

Participants are welcomed to \emph{consult} the jury on any questions on the project. 
If they successfully work on the project, they can get interesting \emph{extra problems}.

A \textit{team} working on this project may consist of any number of participants. 
For every solution \textit{
written for a user} marked with either~`$+$' or `$+.$' a team gets five `beans' (see recommendations in p.~3, `How to write a proof for a user' of 
\url{https://www.mccme.ru/circles/oim/multicomb.pdf}).
The jury may also award extra beans for beautiful solutions, solutions of hard problems, or solutions typeset in \TeX.
The jury has infinitely many beans.
Every team initially has five beans. 
A team may submit a solution \textit{in oral form} or as \textit{written for a developer} if the team has some beans. 
A team loses a bean with every attempt (successful or not).


\comment

Please notify us if you already know solutions of several problems.
If you confirm your knowledge by presenting some solutions, you will be allowed not to receive plus-marks for the  problems, but to use them in solutions of other problems.

Participants (or teams) can
submit their solutions by a personal communication to Egor Riabov at \url{https://mattermost.turgor.ru}.
Please also send him questions and requests for hints on problems which you are stuck with.
 
Participants (or teams) from Serbia and Croatia may submit their solutions to Prof. Rade \v Zivaljevi\'c at \texttt{rade@mi.sanu.ac.rs}.
 

\smallskip
{\bf How to write a proof for a user}

We give some recommendations on how to write a proof that could be included in a mathematical book or research paper (which is a `reliable reference', cf. \url{https://arxiv.org/pdf/2101.03745.pdf}, p. 2).
These recommendations are by no means complete.
You can learn to write proofs (solutions of problems) by trying to write them and discussing your text with a teacher.

See also \url{https://en.wikipedia.org/wiki/KISS_principle}

\url{http://people.apache.org/~fhanik/kiss.html}

{\it A genius makes his own rules, but a `how to' article is written by one ordinary
mortal for the benefit of another...
Most things that an article such as this one can say have at least one
counterexample in the practice of some natural born genius.
Authors of articles such as this one know that, but in the first approximation
they must ignore it, or nothing would ever get done.}

(P. Halmos, How to talk mathematics.)

(1) Only write sentences that make sense.\footnote{For example, none of the following two sentences makes sense:
\[1+2+\ldots+n = \frac{n(n+1)}{2},\]
\[1+2+\ldots+n = \frac{n(n+1)}{2}\text{\ for a positive integer\ }n,\]
because it is not written for which $n$ the statement is stated.
The following statements do make sense:
\[1+2+\ldots+n = \frac{n(n+1)}{2}\text{\ for \textit{every} positive integer\ }n,\]
\[1+2+\ldots+n = \frac{n(n+1)}{2}\text{\ for \textit{some} positive integer\ }n,\]
\[1+2+\ldots+n = 100 \text{\ for \textit{some} positive integer\ }n.\]
However, the second of them is not interesting and the third is not correct.}

(1a) In particular, long sentence usually does not make sense because
it is unclear which exactly parts of long sentence words `and', `or', `then' are tying together.
So break long sentences into short ones.

(2) Introduce notation and each definition explicitly with `define', `denote', `let', `set', `put'.
For example, the phrase `$a=b+c$' without these words means `the previously defined object $a$ equals to the sum of the previously defined objects $b$ and $c$'.

(5) Do not put any part of your solution in parentheses.
Parentheses do not make clear the logical relation between the phrase in and outside the parentheses.
(Parentheses are used for remarks which are not part of the solution.)


\endcomment

\section{One-dimensional cycles in graphs}\label{s:cgra}

The rigorous definition of a graph is given at the beginning of \S\ref{s:1squ}. 
Before that you can work with graphs at an intuitive level. 
  
Denote by
\begin{itemize} 
    \item $[n]:=\{1,2,\ldots,n\}$; 

    \item $K_n$ the complete graph on the set $[n]$ of vertices;
    
    \item $K_{m,n}$ the complete bipartite graph with parts $[m]$ and $[n]'$ (we denote by $A'$ a copy of $A$).
\end{itemize}

A \textbf{simple cycle $v_1v_2\ldots v_k=(v_1,v_2,\ldots,v_k)$ of length $k$} in a graph is a set $v_1v_2, v_2v_3, \ldots, v_kv_1$ of edges such that the vertices $v_1, \ldots, v_k$ are pairwise distinct.
A \emph{simple cycle} is a simple cycle of some length.
We omit `simple' if this word is clear from the context. 
Clearly, any simple cycle is a $1$-cycle.
The definition of a $1$-cycle is given in the introduction. 

The \textbf{sum} (the sum modulo $2$, or the symmetric difference) of sets $A,B$ is
\[
   A+B := (A \cup B) \setminus (A \cap B).
\]

\begin{pr}\label{l:gencom1}
\begin{enumerate}
    \item The sum of $1$-cycles is a $1$-cycle.             

    \item\label{en:gencom1:sum} Any $1$-cycle in $K_n$ is a sum of some cycles of length~$3$.

    \item\label{en:gencom1:tre} If every edge of a $1$-cycle in $K_n$ contains the vertex $n$, then the $1$-cycle is empty. 

    \item\label{en:gencom1:nmb} How many $1$-cycles in $K_n$ are there?

    \item\label{en:gencom1:sumtr} For any vertices $a,b,c,d$ in $K_n$ we have $abc+abd+acd+bcd=0$ (we denote $0:=\emptyset$).

    \item\label{en:gencom1:rel} Any linear relation on cycles of length~$3$ in $K_n$ is a sum of some relations from~(\ref{en:gencom1:sumtr}). 
\end{enumerate}

    E.g. 
    \begin{itemize}
        \item the relation $123+124+134+235+245+345=0$ is the sum of relations
        
        $123+124+134+234=0$ and $234+235+245+345=0$;
    
        \item the sum of boundaries of triangles from any triangulation of the sphere or of the torus (Figure \ref{f:k7-on-torus}) is zero.
    \end{itemize}

\begin{figure}[h]\centering
\includegraphics[width=3cm]{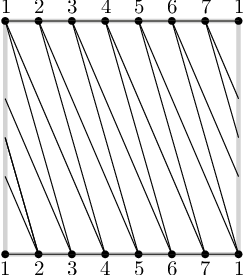} 
\caption{A $7$-vertices triangulation of the torus}\label{f:k7-on-torus}
\end{figure}
 
Rigorously, a \emph{linear relation} is a set of cycles of length~$3$ in $K_n$ such that every edge of $K_n$ is contained in an even number of cycles from this set.    
Thus a rigorous formulation of (\ref{en:gencom1:rel}) is as follows (we identify cycles of length~3 in $K_n$ with 3-element subsets of $[n]$). 
Suppose that $R$ is a set of $3$-element subsets of $[n]$ such that any $2$-element subset of $[n]$ is contained in an even number of those $3$-element subsets that are elements of $R$. 
Then there is a family $R'$ of $4$-element subsets of $[n]$ such that a $3$-element subset $B$ of $[n]$ is contained in $R$ if and only if $B$ is contained in an odd number of those $4$-element subsets that are elements of $R'$.  
\end{pr}

\begin{pr}\label{pr:chl}
\begin{enumerate}
    \item\label{en:chl:simple}
    Any $1$-cycle in a graph is a sum of some simple cycles.
    
    \item\label{en:chl:nochord}
    Any $1$-cycle in a graph is a sum of some simple cycles for which there are no chords.
    A \emph{chord} in a simple cycle is an edge of the graph which does not belong to the cycle (i.e., an edge between two non-consecutive vertices).

    \item\label{en:chl:x-4} Any $1$-cycle in $K_n$ is a sum of some of the following cycles: $123$ and cycles of length~$4$.

    \item\label{en:chl:xneq4} The cycle $123$ is not a sum of some cycles of length~$4$ in $K_n$.
\end{enumerate}
\end{pr}

\begin{pr}\label{p:knn}
\begin{enumerate}
     \item\label{en:knn:sum} Any $1$-cycle in $K_{n,n}$ is a sum of some cycles of length~$4$.

    \item\label{en:knn:tre} If every edge of a $1$-cycle in $K_{n,n}$ contains at least one of the vertices $n$ or $n'$, then the $1$-cycle is empty.

    \item\label{en:knn:numcyc} How many $1$-cycles in $K_{n,n}$ are there?

    \item\label{en:knn:sumrect} For any pairwise distinct $a,b,c\in[n]$ and distinct $u',v'\in[n]'$ we have 
    $ au'bv' + bu'cv' + cu'av' =0$.  

    \item\label{en:knn:rel} Any linear relation on cycles of length~$4$ in $K_{n,n}$ is a sum of some relations from~\eqref{en:knn:sumrect}, and analogous relations $a'ub'v + b'uc'v + c'ua'v = 0$.
    (Rigorously, a \emph{linear relation} is a set of cycles of length~$4$ in $K_{n,n}$ such that every edge of $K_{n,n}$ is contained in an even number of cycles from this set.)
\end{enumerate}
\end{pr}

\begin{pr}\label{pr:orinumber}
    How many $1$-cycles are there
\begin{enumerate}
   \item\label{en:orinumber:tree} in a tree?

   \item\label{en:orinumber:conngraph} in a connected graph with~$V$ vertices and~$E$ edges?
\end{enumerate}
\end{pr}

Denote by $\t{K_n}$ the graph obtained from $K_{n,n}$ by deleting all `diagonal' edges $jj'$, $j\in[n]$. 
E.g., $\t{K_3}$ is the cycle $12'31'23'$ of length $6$. 
This graph naturally appears in Problem \ref{t:2cyc-bij}. 
 
\begin{pr}\label{p:tkn}
\begin{enumerate}
    \item\label{en:tkn:4} If $n\ge4$, then any $1$-cycle in $\t{K_n}$ is a sum of some cycles of length $4$.

    \item\label{en:tkn:6} 
    Any $1$-cycle in $\t{K_n}$ is a sum of some of the cycles $C_{ij'}:=12'31'ij'$ for edges $ij'$ of $\t{K_n}$ such that $i,j>1$ and $(i,j)\ne(3,2)$. 
 
    \item\label{en:tkn:un} The representation of~\eqref{en:tkn:6} is unique.

\end{enumerate}
\end{pr}

A set $B$ of $1$-cycles is called a \emph{base} of a set $A$ of $1$-cycles, if every $1$-cycle from $A$ has a unique representation as a sum of some $1$-cycles from $B$. 

\begin{pr}[symmetric $1$-cycles]\label{p:symgra}*
Denote by $t\colon\t{K_n}\to\t{K_n}$ the symmetry (involution) switching the parts, i.e., switching $j$ and $j'$ 
for every $j\in[n]$.  
Denote by $tQ$ the $1$-cycle symmetric to a $1$-cycle $Q$.   
E.g. $t\t{K_3}=\t{K_3}$ and $t(12'31'ij')=1'23'1i'j$.   
A $1$-cycle $Q$ is called \emph{$t$-symmetric} if $tQ=Q$.
E.g. $\t{K_3}$ is $t$-symmetric.

    (a) How many $t$-symmetric $1$-cycles are there in $\t{K_n}$?
    
    (b) Any $t$-symmetric $1$-cycle in $\t{K_n}$ is the sum of some of the following $1$-cycles: $\t{K_3}$ and $Q+tQ$ for cycles $Q$ of length $4$.
    
    (c) The set of all $1$-cycles in $\t{K_n}$ has a base consisting of (the $t$-symmetric $1$-cycle) $\t{K_3}$ and pairs of $1$-cycles going to each other under $t$.  
\end{pr}


\begin{pr}[integer $1$-cycles]\label{3-homolzcx}*
    Let $K$ be a graph with oriented edges.
    An assignment of integers to oriented edges of $K$ is a (simplicial) \emph{integer $1$-cycle} if for every vertex the sum of integers assigned to incoming edges equals the sum of integers assigned to outcoming edges (Kirchhoff rule).

    (a,b)~State and prove the analogues of Problems~\ref{pr:orinumber}.\ref{en:orinumber:tree},\ref{en:orinumber:conngraph} for integer $1$-cycles.
    
    (c)~For given set $\omega$ of orientations on edges of $K$ denote by $H_1^\omega(K;\Z)$ the set of the integer $1$-cycles, with the componentwise sum operation.
    For different sets $\omega,\omega'$ of orientations there is a 1--1 correspondence $\varphi\colon H_1^\omega(K;\Z)\to H_1^{\omega'}(K;\Z)$ such that $\varphi(x+y)=\varphi(x)+\varphi(y)$ for every $x,y\in H_1^\omega(K;\Z)$ (i.e., the \emph{groups} $H_1(K;\Z)$ for different sets of orientations are \emph{isomorphic}).
\end{pr}

\begin{pr}[cohomology]\label{p:coh}*
Fix an assignment of $+$ or $-$ to edges of a graph.
For any vertex one can inverse the sign of every edge containing this vertex.
    
(a) For a tree, by the above operations for different vertices one can obtain any assignment from any other assignment.


(b) Consider a connected graph with $V$ vertices and $E$ edges.
    Find a maximal number of assignments of $+$ or $-$ to edges such that none of these assignments can be obtained from any other by the above operations.
\end{pr}


\section{One-dimensional cycles in the square of a graph}\label{s:1squ}

A \textbf{graph} is a pair $(V,E)$ of a finite set~$V$ and a set~$E$ of $2$-element subsets (called `edges') of~$V$.

\emph{In this text $K$ is a graph. We sometimes denote edge $\{a,b\}$ shortly by $ab$.}

Intuitively, a \textit{configuration} is a pair of ants (one red and one blue) sitting at two vertices of a graph (ants may occupy the same vertex). 
Configurations are called \textit{adjacent} if one can be obtained from the other by one of the ants moving along one edge. 
This is rigorized as follows. 
Vertices of the graph $K^{\square 2}$ are ordered pairs $(a,b)$ of vertices of the graph $K$.
If vertices $b$ and $c$ of $K$ are joined by an edge, then the vertices $(a,b)$ and $(a,c)$ of $K^{\square 2}$ are joined by an edge denoted by $(a,bc)$, and vertices $(b,a)$ and $(c,a)$ are joined by an edge denoted by $(bc,a)$.
There are no other edges in $K^{\square 2}$.
 
E.g., the following are simple cycles in $K^{\square 2}$ (recall from \S\ref{s:cgra} that a simple cycle is denoted by vertices it passes through): 

$\bullet$ the \emph{boundary} 
$$ab\Box uv := (a,u)(b,u)(b,v)(a,v),$$ 
for edges $ab$ and $uv$ in $K$
(the product $ab\times uv$ is a rectangle, see \S\ref{s:prod}, and $ab\Box uv=\partial(ab\times uv)$ is its boundary);

$\bullet$ the \emph{diagonal}, \emph{off-diagonal} and \emph{antidiagonal cycles} 
$$\diag(123):=(1,1)(1,2)(2,2)(2,3)(3,3)(3,1),\quad(1,2)(1,3)(2,3)(2,1)(3,1)(3,2),$$ 
$$\text{and}\quad(1,1)(2,1)(2,3)(3,3)(3,2)(1,2)\quad\text{in}\quad K_3^{\square 2};$$ 

$\bullet$  the \emph{left}, \emph{right}, \emph{diagonal}, \emph{off-diagonal}, and \emph{antidiagonal cycles} 
$$a\times C:=(a,v_1)\ldots(a,v_k),\quad C\times a:=(v_1,a)\ldots(v_k,a),$$
$$\diag C:=(v_1,v_1)(v_1,v_2)(v_2,v_2)\ldots(v_k,v_k)(v_k,v_1),\quad (v_1,v_2)(v_1,v_3)(v_2,v_3)\ldots(v_k,v_1)(v_k,v_2),$$ 
$$\text{and}\quad(v_1,v_1)(v_2,v_1)(v_2,v_k)\ldots(v_k,v_2)(v_1,v_2)$$  
for a vertex $a$ and a simple cycle $C=v_1\ldots v_k$ in $K$; 


$\bullet$ the \emph{triodic cycle} $(1,3)(1,1')(1,2)(1',2)(3,2)(3,1')\ldots$ in $K_{3,1}^{\square 2}$, where dots denote the part symmetric to the written part; a \emph{triodic cycle} is also the analogous cycle corresponding to a $K_{3,1}$-subgraph of $K$.  

We consider the symmetry (involution) of $K^{\square 2}$ switching the factors (i.e., switching points $(x,y)$ and $(y,x)$), and the corresponding map on $1$-cycles.   
(For an application of symmetric 1-cycles 
see e.g. 
\cite[\S1.6]{Sk}.)

\begin{pr}[riddle]\label{hatg2} Find graphs $K^{\square 2}$ for \quad 

(a) $K=K_{2,1}$ a path on three vertices; \quad 
(b) $K=K_3$ a cycle on three vertices; \quad 

(c) $K=K_{3,1}$ a triod; \quad (d) $K=K_4$.
\end{pr}

\begin{pr}\label{p:treek2}
(a) Any 1-cycle in $K_{2,1}^{\square 2}$ is a sum of some boundaries. 

(b) The triodic cycle in $K_{3,1}^{\square 2}$ is a sum of some boundaries. 

(c) Any 1-cycle in $K_{3,1}^{\square 2}$ is a sum of some boundaries. 

(d) If $K$ is a tree, then any 1-cycle in $K^{\square 2}$ is a sum of some boundaries. 

{\it Hint.} 
Prove that \emph{if $K$ and $L$ are trees, then any $1$-cycle in the graph $K\Box L$} (defined below) \emph{is a sum of some boundaries}. 
Prove this by induction, using deletion of a leaf vertex. 

(e) Is some 
left cycle a sum of some boundaries?  

(f) Is some 
diagonal cycle a sum of some boundaries?  
\end{pr}

The following notions are useful, in particular, for Problems \ref{p:treek2}.def. 

Let $K$ and $L$ be graphs. 
Vertices of the \emph{product graph} $K\Box L$ are ordered pairs $(a,b)$ of vertices $a$ of $K$ and $b$ of $L$.
If vertices $b$ and $c$ of $L$ are joined by an edge, then the vertices $(a,b)$ and $(a,c)$ of $K\Box L$ are joined by an edge denoted by $(a,bc)$. 
If vertices $b$ and $c$ of $K$ are joined by an edge, then the vertices $(b,a)$ and $(c,a)$ of $K\Box L$ are joined by an edge denoted by $(bc,a)$.  
There are no other edges in $K\Box L$.

\emph{The left projection} $C_y$ of a 1-cycle $C$ in $K\Box L$ is the set of all edges~$\sigma$ in~$L$ such that there is an odd number of vertices~$a$ in~$K$ such that $a\times\sigma\in C$. 
(If $C$ is a cycle, then $C_y$ is the set of all edges passed by the blue ant an odd number of times.) 
\emph{The right projection} $C_x$ is defined analogously.  

\begin{pr}\label{p:diag}  
(a) The diagonal cycle in $K_3^{\square 2}$ is a sum of a left cycle, a right cycle, and some boundaries. 

(b) The off-diagonal cycle in $K_3^{\square 2}$ is a sum of the diagonal cycle and some boundaries. 

(c) The antidiagonal cycle in $K_3^{\square 2}$ is a sum of the diagonal cycle and some boundaries. 

(d) The \emph{symmetrized cycle} $a\times C+C\times a$ is the sum of $\diag C$ and some boundaries.  

(e) For $K=K_4$ is the left cycle $1\times234$ a sum of some diagonal cycles and some boundaries? 
\end{pr}

\begin{pr}\label{p:k21} Find the number of $1$-cycles up to boundaries in

(a)  $K_3^{\square 2}$; \quad (b)  $K_{2,2}^{\square 2}$; \quad 
(c)  $K_{2,3}^{\square 2}$; \quad (d)  $K_4^{\square 2}$. 

{\it Hint.} 
The number of $1$-cycles up to boundaries is the maximal number of $1$-cycles such that none of them is the sum of any other and some boundaries. 
For some parts of this problem the K\"unneth theorem \ref{stcycles2}.c will be useful. 
\end{pr}

Two 1-cycles $C,C'$ in $K^{\square 2}$ are said to be \emph{homologous} (or congruent modulo boundaries) if $C+C'$ is a sum of some boundaries. 
Notation: $C\sim C'$. 

\begin{pr}\label{stcycles2}
(a)  If $K$ is connected and $T$ is a tree, then for any $1$-cycle $C$ in $K\Box T$
and vertex $a$ in $T$ there is a unique 1-cycle $C_1$ in $K$ such that $C\sim C_1\times a$. 
Moreover, $C\sim C_x\times a$. 

(b) If $Z$ is a 1-cycle in $K^{\square 2}$ such that $Z_x=0$, then $Z\sim a\times Z_y$ for any vertex $a$.

(c) (K\"unneth theorem) If $K$ is connected, then for any $1$-cycle $C$ in $K^{\square 2}$ and vertex $a$ in $K$ there are unique 1-cycles $C_1,C_2$ in $K$ such that $C\sim C_1\times a+a\times C_2$. 
Moreover, $C\sim C_x\times a + a\times C_y$.
\end{pr}

See more information on $K\Box L$ in \cite{CPG}. 

The following (and $K^{\underline2}$ defined in \S\ref{s:cydepr}) is a graph-theoretical analogue of the set of arrangements.
Consider two ants on $K$ as above, but now forbidden to sit in the same vertex. 
Rigoriously speaking, vertices of the graph $K^{\square \underline2}$ are ordered pairs of distinct vertices of~$K$.
Vertices of $K^{\square \underline2}$ are joined by an edge in $K^{\square \underline2}$ if they are joined by an edge in $K^{\square 2}$.

\begin{pr}\label{hatg} 
(riddle) Find graphs $K^{\square \underline2}$ in cases (a,b,c,d) of Problem~\ref{hatg2}.
\end{pr}

\begin{pr}\label{p:knn21} Find the number of $1$-cycles in 

(a)  $K_3^{\square \underline2}$; \quad (b)  $K_{2,2}^{\square \underline2}$; \quad 
(c)  $K_{2,3}^{\square \underline2}$; \quad (d)  $K_4^{\square \underline2}$; \quad (e) $K_{3,3}^{\square \underline2}$; \quad (f) $K_5^{\square \underline2}$. 

(a'-f') Find the number of $1$-cycles in $K^{\square\underline2}$ up to boundaries (contained) in $K^{\square \underline2}$ (i.e., up to boundaries corresponding to pairs of non-adjacent
edges).  
 
{\it Hints.} (a'-d') If a set of boundaries in $K^{\square\underline2}$ has zero sum, then the set is empty. 
 
(e'), (f') The sum of all boundaries in $K^{\square\underline2}$ is zero. 
This is the only non-empty set of boundaries in $K^{\square\underline2}$ that has zero sum. 

\end{pr}
 
   
\begin{pr}\label{p:stcycles}* 
(a) If a sum of boundaries is symmetric, it is a \emph{symmetrized boundary}, i.e., a sum of $B+B'$ for some pairs of boundaries $B,B'$ symmetric to each other.

(b) 
Any symmetric $1$-cycle in $K^{\square 2}$ is a sum of some symmetrized cycles and some symmetrized boundaries. 
(Hence it is a sum of some diagonal cycles and some boundaries; note that a diagonal cycle is  symmetric up to adding boundaries.)  
 
{\it Hint.} 
Follows by (a) and the K\"unneth Theorem \ref{stcycles2}.c. 

(c) 
For $K$ connected, the map $C\mapsto C\times a + a\times C$ defines a 1--1 correspondence between $1$-cycles in $K$ and symmetric $1$-cycles in $K^{\square 2}$ up to symmetrized boundaries. 
The inverse correspondence is given by $C\mapsto C_x=C_y$. 
 
(d) For $K$ connected, the number of symmetric $1$-cycles in $K^{\square 2}$ is uniquely defined by the number $V$ of vertices and the number $E$ of edges. \cite[Propoistion 5.3]{DMN+} 
\end{pr}

\begin{pr}\label{stcyclesdel}* Let $K$ be a connected graph. 

(a) A triodic cycle is not a sum of some boundaries in $K_{3,1}^{\square \underline2}$. 

(b) Is an antidiagonal cycle a sum of a 1-cycle in $K^{\square \underline2}$ and some boundaries (in $K^{\square 2}$)? 

(c) (conjecture) In $K^{\square\underline2}$ any off-diagonal cycle is not a sum of symmetrized cycles, and boundaries. 

(d) In $K_4^{\square \underline2}$ is any off-diagonal cycle a sum of symmetrized cycles, triodic cycles, and boundaries? 

(e) In $K^{\square\underline2}$ is any symmetrized cycle a sum of off-diagonal cycles, triodic cycles, and boundaries.  


(f) (conjecture) In $K^{\square\underline2}$ any 1-cycle is a sum of left cycles, right cycles, off-diagonal cycles, triodic cycles, and boundaries. 

(g) In $K^{\square \underline2}$ any symmetric 1-cycle  is a sum of off-diagonal cycles (redefined in a symmetric way), triodic cycles, and boundaries. \cite{Bo}  
\end{pr}

\section{Cycles in hypergraphs}\label{s:chygra}
 
Higher-dimensional cycles in hypergraphs appear, in particular, as relations on 1-cycles in a graph, see Assertions \ref{l:gencom1}.\ref{en:gencom1:sumtr}\ref{en:gencom1:rel}.
 
\begin{pr}[cf. Problem \ref{l:gencom1}]\label{l:gencom}
A \textit{$2$-cycle} is a set $C$ of $3$-element subsets (called `faces') of $[n]$ such that every $2$-element subsets of $[n]$ is contained in an even number of subsets from $C$. 
E.g. the empty set is a $2$-cycle. 
For a $4$-element subset $A\subset[n]$ let the \textit{tetrahedron} $T_A$ be the set of all 3-element subsets of $A$. 
In other words, for pairwise distinct $a,b,c,d\in [n]$ define the \textit{tetrahedron} 
$$T_{\{a,b,c,d\}}:=\bigl\{\{a,b,c\},\{a,b,d\},\{a,c,d\},\{b,c,d\}\bigr\}.$$ 
Clearly, any tetrahedron is a $2$-cycle. 

\begin{enumerate}
    \item The sum of $2$-cycles is a $2$-cycle. 

    \item\label{en:gencom:tre} If every face of a $2$-cycle contains the number $n$, then the $2$-cycle is empty. 

    \item\label{en:gencom:gen} Any $2$-cycle is a sum of some tetrahedra. 

    \item How many $2$-cycles for $[n]$ are there?
     
    \item\label{en:gencom:sumtetra} For any $5$-element subset $A\subset[n]$ 
we have $\sum\limits_{j\in A}T_{A-\{j\}}=0$.

    \item\label{en:gencom:rel} Any linear relation on tetrahedra is a sum of some relations from~\eqref{en:gencom:sumtetra}. 
    (Give a rigorous formulation analogous to Assertion \ref{l:gencom1}.\ref{en:gencom1:rel}.) 
\end{enumerate}

\end{pr}

\begin{pr}\label{pr:rook} 
Recall that $[n]^\ell$ is the set of all vectors of length $\ell$ whose entries are numbers from $[n]$.  
A~\emph{line} is a subset of $[n]^\ell$ given by fixing all coordinates except one. 
A \textit{rook cycle} is a~subset of~$[n]^\ell$ containing an even number of vertices of every line. 
E.g. the empty set is a rook cycle. 
A \textit{parallelepiped} is a subset $P_1\times\ldots\times P_{\ell}\subset[n]^\ell$, where $P_i$ are $2$-element subsets of $[n]$. 
Clearly, any parallelepiped is a rook cycle.  

(a) The sum of rook cycles is a rook cycle. 

(b) If a rook cycle contains no elements from $[n-1]^\ell$, then the rook cycle is empty.

(c) Any rook cycle is a sum of some parallelepipeds.

(d) How many rook cycles for $[n]^\ell$ are there?

(e) For pairwise distinct $a,b,c\in[n]$ and parallelepiped $P\subset[n]^{\ell-1}$ we have 
\linebreak
$P\times\{a,b\}+P\times\{b,c\}+P\times\{c,a\}=0$.  

(f) Any linear relation on parallelepipeds in $[n]^\ell$ is a sum of some relations from (e), and relations obtained from that by permutations of coordinates.
(A rigorous formulation is analogous to Assertions~\ref{p:knn}.\ref{en:knn:rel} and \ref{l:gencom}.\ref{en:gencom:rel}.)

{\it Hint.}  First try the case $\ell=2$, cf. Problem \ref{p:knn}. 
\end{pr}

\begin{pr}[riddle]\label{p:high} 
    Invent and prove a higher-dimensional analogue of Problem \ref{l:gencom}. 
\end{pr}

\begin{pr}\label{p:cychyp} * 
A \emph{$2$-hypergraph} (a two-dimensional hypergraph, or a $3$-uniform hypergraph) is a pair $(V,F)$ of a finite set~$V$ and a set~$F$ of $3$-element subsets (`faces') of~$V$.
A (simplicial) \emph{$2$-cycle} (modulo $2$) in a $2$-hypergraph $(V,F)$ is a set $C$ of faces  such that every $2$-element subset of $V$ is contained in an even number of subsets from $C$. 
    
(a)   There are two \emph{connected} hypergraphs having the same numbers of vertices, edges (i.e., $2$-element subsets of faces), and faces, but different numbers of $2$-cycles. (Cf. Problem \ref{pr:orinumber}.\ref{en:orinumber:conngraph}.) 

(b) Let $L=(V,F)$ be a 
hypergraph. 
Denote by $E$ the set of edges, by $b_0$ the number of connected components, 
by $2^{b_1}$ the number of $1$-cycles up to sums of some boundaries of faces, 
by $2^{b_2}$ the number of $2$-cycles.
Then $b_0-b_1+b_2=|V|-|E|+|F|$. 
    
\end{pr}


\section{Geometric digression: Cartesian products of graphs}\label{s:prod}
 
First we give an intuitive definition of a cylinder, and then a rigorous one. 

Take rectangle bands $aa'bb'$ correponding to edges $ab$ of a graph $K$.
Glue the end segments of the bands corresponding to the same vertex, so that the hatched letters would stick together with the hatched ones.
The obtained two-dimensional figure is called a \emph{cylinder} over graph $K$.
 
\begin{figure}[h]
\centerline{\includegraphics[width=3cm]{real-83.mps}\qquad
\includegraphics[width=3.5cm]{real-23.mps}\qquad
\includegraphics[width=3cm]{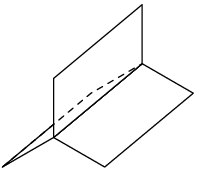}\qquad 
\includegraphics[width=5cm]{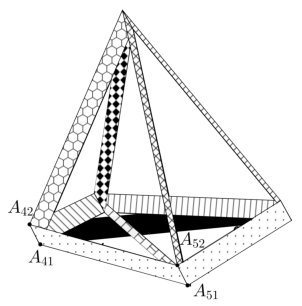}
}
\caption{The cylinders over $K_2$, $K_3$, $K_{3,1}$, $K_5$} 
\label{f-mnre}
\end{figure}

Recall that $\R^d$ is the $d$-dimensional Euclidean space (for $d=2$ and $d=3$ this is the usual plane and space which one studies in elementary geometry class). 
The \emph{cylinder} over a subset $U\subset\R^d$ is 
$$U\times K_2:=\{(x,t)\in\R^{d+1}\ :\ x\in U,\ t\in[-1,1]\}.$$
For instance, the cylinders over $K_2$, $K_3$, $K_{3,1}$, and $K_5$ are shown in Figure~\ref{f-mnre}; 
the cylinder over $K_{k,1}$ `looks like' the book with $k$ sheets. 

\begin{pr}\label{p:realbook} 
Any graph can be drawn without edges crossing on a book with a certain number of sheets (Figure~\ref{f-mnre}, $K_{3,1}$) depending on the graph. 
More precisely, for any~$n$ there exists an integer~$k$, as well as $n$~points and $n(n-1)/2$~non-self-intersecting polygonal lines on a book with $k$~sheets such that 

$\bullet$ every pair of points is joined by a polygonal line, and 

$\bullet$ no polygonal line intersects the interior of another polygonal line.

\end{pr}

\begin{pr}[riddle]\label{p:realcyl} 
The cylinder over any graph is realizable without self-intersections in $\R^3$. 
\end{pr}

The \emph{(geometric) product} of subsets $U,V\subset\R^d$ is 
$$U\times V:=\{(x,y)\in\R^{2d}\ :\ x\in U,\ y\in V\}.$$
In particular, the \emph{square} of a subset $U\subset\R^d$ is $U^2=U\times U$.  
For instance, the square $K_2^2$ is the ordinary square in the plane, and the square $K_3^2$ (called the \emph{torus}) is shown in Figure~\ref{f:k5i}, middle.

\begin{pr}[riddle]\label{p:cyl} Draw in $\R^3$ 

(a) the product $K_3\times K_{3,1}$ without self-intersections. 

(b)  the square $K_{3,1}^2$; self-intersections are allowed. 
\end{pr}

\textbf{Remark.}
Sometimes the word `graph' is used for `the body of a graph' defined as follows. 
Suppose that a subset of $\R^3$ is in 1--1 correspondence with the set of vertices of $K$, and no segment corresponding to an edge of $K$ intersects the interior of any other such segment. 
The \textit{body} of $K$ is the union of the subset, and of all such segments. 

In this language $K^{\square 2}$ is the set (the configuration space) of ordered pairs $(x,y)$ of points of a graph $K$ (or rather of its body) such that either $x$ or $y$ is a vertex of $K$. 
The set $K^{\square 2}$ is a union of a finite number of segments, i.e., is a graph.

Analogously, $K^{\square \underline2}$ is the set (the configuration space) of ordered pairs $(x,y)$ of points of a graph $K$ (or rather of its body) such that either $x$ or $y$ is a vertex of $K$, and $x,y$ do not belong to the same edge.
The set $K^{\square \underline2}$ is a union of a finite number of segments, i.e., is a graph.


\section{Two-dimensional cycles in the square of a graph}\label{s:cydepr}

Let \emph{the (combinatorial) product} of graphs $K$ and $L$ be 
$$K\times L := \{(\sigma,\tau)\ :\ \sigma\text{ is an edge of }K,\ \tau\text{ is an edge of }L\}.$$
Cells $(\sigma, \tau)$ and $(\gamma, \delta)$ in $K\times L$ are said to be \emph{adjacent} if either 
    
$\bullet$ $\sigma=\gamma$, and $\tau, \delta$ have a common vertex, or 

$\bullet$ $\tau=\delta$, and $\sigma, \gamma$ have a common vertex.  

(i.e., if the cells `share a common edge' in $K\Box L$).
If we glue adjacent cells along their `common edges' in $K\Box L$, we obtain the geometric product discussed in \S\ref{s:prod}. 

Let \emph{the square} of $K$ be $K^2 := K\times K$. 

A \textbf{cellular $2$-cycle} (in $K^2$) is a subset $C\subset K^2$ such that for every vertex $a$ and edge $\beta$ of $K$

$\bullet$  there is an even number of edges $\alpha$ of $K$ such that $\alpha\ni a$ and $(\alpha,\beta)\in C$, and 

$\bullet$  there is an even number of edges $\alpha$ of $K$ such that $\alpha\ni a$ and $(\beta,\alpha)\in C$. 

In other words, $C\subset K^2$ is a cellular 2-cycle if every edge in $K^{\square 2}$ is `contained' in an even number of cells $(\alpha,\beta)\in C$.  

\begin{pr}\label{p:ce2cy} 
(a) For subset $C$ of $K^2$, the sum of boundaries $\sigma\Box\tau$ over $(\sigma,\tau)\in C$ is zero if and only if $C$ is a cellular $2$-cycle.  

(b) A subset $C\subset K^2$ is a cellular $2$-cycle if and only if for every edge $\sigma$ of $K$ both sets (`column section' and `row section')
\[
    C_{\sigma,\cdot}:=\bigl\{\tau\ :\ (\sigma,\tau)\in C\bigr\}\quad\text{and}\quad 
    C_{\cdot,\sigma}:=\bigl\{\tau\ :\ (\tau,\sigma)\in C\bigr\}
\]
are $1$-cycles in $K$. 
\end{pr}

\begin{pr}\label{t:torus} 
(a) The product of any two simple cycles in $K$ (the \emph{torus}) is a cellular $2$-cycle. 

(b) The sum modulo $2$ of cellular $2$-cycles is a cellular $2$-cycle. 

(c) Does there exist a tree $T$, and a non-empty cellular $2$-cycle in $T^2$? 


(d) Does there exist a graph $K$, a tree $T\subset K$, and a non-empty cellular $2$-cycle contained in $\overline T:=T\times K\cup K\times T$? 
\end{pr} 

\begin{pr}\label{p:k2}  
Find the number of cellular $2$-cycles in 

(a)  $K_3^2$; \quad (b)  $K_{2,2}^2$; \quad 
(c)  $K_{2,3}^2$; \quad (d)  $K_4^2$.     

{\it Hint.} For some parts of this problem the following problems will be useful. 
\end{pr}

Intuitively (and not precisely) speaking, the \emph{combinatorial deleted product} $K^{\underline2}$ of $K$ is a set of ordered pairs $(x,y)$ of points in the graph $K$ such that $x,y$ do not belong to adjacent edges.
The set $K^{\underline2}$ could be represented as a union of rectangles.
Rigorously, let the \emph{combinatorial deleted product} of $K$ be     
$$K^{\underline2} := \{(\sigma,\tau)\ :\ \sigma,\tau\text{ are non-adjacent
edges of }K\}.$$ 

{\bf Remark.} This set and its generalizations have many applications, notably to algorithm recognizing realizability of hypergraphs in higher-dimensional Euclidean spaces, see surveys \cite{Sk06, Sk18}. 
The subset $K^{\underline2}\subset K^2$ is smaller than $K^2$, and so is sometimes easier-to-draw. 
E.g. $K_2^{\underline2}=K_3^{\underline2}=K_{n,1}^{\underline2}=\emptyset$. 


Glue adjacent cells along their `common edges' 
in $K^{\square 2}$. 
Then the combinatorial deleted product of the path on 5 vertices, of the cycle on 5 vertices, of $K_4$, of $K_{3,3}$, and of $K_5$ look like the disjoint union of two disks, the annulus, the cuboctahedron without triangular faces (Figure \ref{f:del3}), the sphere with 4 handles, and the sphere with 6 handles, respectively. 

(For $K=K_{3,3}$ and $K=K_5$ one proves that $K^{\underline 2}$ looks like a connected orientable two-dimensional surface, and uses the classification theorem for two-dimensional surfaces; see details in \cite[the text after 3.4.1]{Sa91}.)

\begin{pr}\label{t:dpcycle} 
(a) Are there any non-empty cellular $2$-cycles in $K^{\underline2}$ for a \emph{cycle} $K$? 

(b) Analogous question for a \emph{wheel} $K$, i.e., for the graph with the set $\{0\}\cup[n]$ of vertices, and edges $\{n,1\}$, $\{0,j\}$ and $\{j,j+1\}$ for $j\in[n-1]$. 

(c) The subset $K_{3,3}^{\underline2}$ of $K_{3,3}^2$ is a cellular $2$-cycle. 

(d) Analogous statement for $K_5$. 

(e) There is a graph $K$ such that $K^{\underline2}$ has a non-empty proper subset that is a cellular $2$-cycle. 
\end{pr} 

\begin{pr}\label{t:kunn} 
(a) Express $K_{3,3}^{\underline2}$ as a sum of some tori. 

(b) Express $K_5^{\underline2}$ as a sum of some tori. 

(c) Any cellular $2$-cycle (in $K^2$) is a sum of some tori. 

(d) How many cellular $2$-cycles in $K^2$ are there, for a connected graph $K$ with $V$ vertices and $E$ edges?

{\it Hint to (c, d).} The $1$-cycle $\widehat \varphi \sigma$ is defined in the
\simple{solution}\sophi{\hyperref[en:sol:orinumber:conngraph]{solution}}
of Problem \ref{pr:orinumber}.\ref{en:orinumber:conngraph}. 
By answer to Problem~\ref{t:torus}.d any cellular $2$-cycle $C$ in $K^2$ equals $\displaystyle\sum_{(\sigma,\tau)\in C-\overline T}\widehat \varphi \sigma\times\widehat \varphi \tau$, and is uniquely defined by its cells outside $\overline T$.
\end{pr} 

\emph{Denote by $H_1(K)$ the set of all $1$-cycles in $K$, with the sum operation.} 

\emph{Denote by $H_2(K^2)$ the set of all cellular $2$-cycles in $K^2$, with the sum operation.} 


\begin{pr}[K\"unneth theorem]\label{t:kunneth} 
(a) There is a base $C_1,\ldots,C_q$ in $H_1(K)$ such that $C_i\times C_j$, $i,j\in[q]$, is a base in $H_2(K^2)$. 

(b) If $C_1,\ldots,C_q$ is a base in $H_1(K)$, then $C_i\times C_j$, $i,j\in[q]$, is a base in $H_2(K^2)$. 

(I.e., $H_2(K^2)\cong H_1(K)\otimes H_1(K)$. Base is defined analogously to the definition before Problem~\ref{p:symgra}.)
\end{pr} 

\begin{pr}\label{t:kunnethsym} 
(a) The product of any two vertex-disjoint cycles in $K$ (i.e., two cycles which do not have common vertices) is a cellular $2$-cycle in $K^{\underline2}$. 

(b) Is any cellular $2$-cycle in $K^{\underline2}$ a sum of some products of vertex-disjoint cycles?
\end{pr}

\begin{pr}\label{t:2cyc-bij}
\begin{enumerate}
    \item\label{en:2cyc-bij:num} How many cellular $2$-cycles in $K_{n,n}^{\underline2}$ are there?

    \item\label{en:2cyc-bij:cyc} There is a 1--1 correspondence $H_2(K_{n,n}^{\underline2})\to H_2(\t{K_n}^2)$.
 
    \item\label{en:2cyc-bij:compl} 
    There is a 1--1 correspondence between $K_{n,n}^{\underline2}$ and $\t{K_n}^2$ that respects the adjacence.
\end{enumerate}

{\it Hint:} solve parts (a-c) in the reverse order.
\end{pr}
 
Consider the symmetry (involution) of $K^2$ switching the factors (i.e., switching points $(x,y)$ and $(y,x)$), and the map induced by this symmetry on cellular $2$-cycles.   

\begin{pr}\label{p:tsym} *  
(a) For any simple cycles  $Q,R$ in $K$ the \emph{symmetrized torus} $Q\times R+R\times Q$ is a symmetric cellular $2$-cycle. 

(b) The sum modulo 2 of symmetric cellular $2$-cycles is a symmetric cellular $2$-cycle. 

(c) Is any symmetric cellular $2$-cycle a sum of some symmetrized tori?  
 
(e) For $K$ connected, the number of symmetric cellular $2$-cycles in $K^2$ uniquely defined by the number $V$ of vertices and the number $E$ of edges. \cite[Proposition 9.2]{DMN+} 
\end{pr} 

\begin{pr}\label{l:h2sym} * 
(a) Any symmetric cellular $2$-cycle contained in $K^{\underline2}$ is a sum of some symmetrized tori. 

(b) How many symmetric cellular $2$-cycles in $K_{n,n}^{\underline2}$ are there?

(c) The 1--1 correspondence of Assertion \ref{t:2cyc-bij}.\ref{en:2cyc-bij:cyc} moves cellular $2$-cycles in $K_{n,n}^{\underline2}$ transposed by the symmetry to cellular $2$-cycles in $\t{K_n}^2$ transposed by the symmetry $t^2$.  
Here the symmetry $t:\t{K_n}\to\t{K_n}$ is defined in Problem \ref{p:symgra}, and $t^2(x,y):=(tx,ty)$. 

(d) Any symmetric cellular $2$-cycle in $K_{n,n}^{\underline2}$ is a sum of some of the following ones: 

$\bullet$ symmetrized tori $Q\times R+R\times Q$ for vertex-disjoint cycles $Q,R$ of length $4$ in $K_{n,n}$; 

$\bullet$ the combinatorial deleted product $K_{3,3}^{\underline2}$ of the subgraph $K_{3,3}$ of $K_{n,n}$. 
 
(e) Any $t^2$-symmetric cellular $2$-cycle in $\t{K_n}^2$ is a sum of some of the following ones: $t^2$-symmetrized tori $Q\times R+tQ\times tR$, and $\t{K_3}^2$. 

(f) Is the analogue of (d) for $K_n$ and the subgraph (homeomorphic to) $K_5$ of $K_n$ correct? 

{\it Hints.} (d) One needs K\"unneth Theorem \ref{t:kunneth}, Assertion~\ref{t:2cyc-bij}.\ref{en:2cyc-bij:cyc} and (e). 

(e) One needs Assertion \ref{p:symgra}.c. 
\end{pr}


\begin{pr}\label{p:sawrong} * 
(a) The following assertion \cite[3.4.2]{Sa91} is wrong, even for connected graphs:  
there is a cellular 2-cycle $C$ in $K^{\underline2}$ such that any cellular 2-cycle in $K^{\underline2}$ is a sum of some of the following cellular 2-cycles: $C$ and products of vertex-disjoint cycles. 

(b) (open problem) Is the assertion of (a) correct for 3-connected graphs? 

(c) Any cellular $2$-cycle in $K^{\underline2}$ is a sum of some products of vertex-disjoint cycles, and some combinatorial deleted products of subgraphs homeomorphic to $K_5$ or to $K_{3,3}$. 

(d) (riddle) How many cellular $2$-cycles in $K^{\underline2}$ are there, for a connected graph $K$ with $V$ vertices and $E$ edges (use additional data of $K$, if required)?

{\it Hints.} (a) Take disjoint union of two copies of $K_5$.  
Make the disjoint union connected by joining the copies by an edge.

\end{pr}


\section{Answers, hints and solutions}\label{s:answ}


\textbf{\ref{l:gencom1}.} 
(\ref{en:gencom1:sum}) The assertion follows from Assertion \ref{pr:chl}.\ref{en:chl:nochord}. 

\textit{Remark.} A specific representation of a $1$-cycle as a sum of some cycles of length~$3$ is given in (\ref{en:gencom1:nmb}). 

(\ref{en:gencom1:nmb}) \emph{Answer:} $2^{\frac{(n-1)(n-2)}{2}}$.

\emph{Hint.} An arbitrary $1$-cycle $C$ in $K_n$ equals the sum $\widehat C := \sum\limits_{ij \in C,\ i,j < n} ijn$ of cycles $ijn$ over all edges $ij\in C$ not containing the vertex $n$.
This follows because every edge in $ C + \widehat C $ contains the vertex $n$, so by~(\ref{en:gencom1:tre}) we have $C + \widehat C = 0$, i.e., $C=\widehat C$.

Moreover, such a representation is unique, see details in the proof of Assertion \ref{pr:orinumber}.\ref{en:orinumber:conngraph}.
 
(\ref{en:gencom1:rel}) In the rigorous formulation this is equivalent to Assertion \ref{l:gencom}.\ref{en:gencom:gen}. 

\smallskip
\textbf{\ref{pr:chl}.} 
(\ref{en:chl:nochord}) Any 1-cycle is a sum of some simple cycles. 
Any chord for a~simple cycle of length~$l$ gives a representation of the cycle as the sum of two simple cycles of lengths less than~$l$.

(\ref{en:chl:x-4}) The sum of two cycles of length~$3$ having a common edge is a cycle of length~$4$. 
So if a cycle $\alpha$ of length~$3$ is a sum of $123$ and cycles of length 4, then any cycle of length~$3$ having a common edge with $\alpha$ is such. 
Therefore, all cycles of length~$3$ are such. 
Now (\ref{en:chl:x-4}) follows by Assertion \ref{l:gencom1}.\ref{en:gencom1:sum}.

(\ref{en:chl:xneq4}) For a finite set $X$ denote by $|X|_2$ the size of $X$ modulo $2$.
Then for any $1$-cycles $C,C'$ we have $|C+C'|_2 = |C|_2 + |C'|_2$.
 
\medskip
\emph{We often shorten $\varphi(x)$ to $\varphi x$ for the image of $x$ under $\varphi$, and $\{\sigma\}$ to $\sigma$ for $1$-element subset.} 

\smallskip
\textbf{\ref{p:knn}.}
(\ref{en:knn:sum}) The assertion follows from Assertion \ref{pr:chl}.\ref{en:chl:nochord}. 

\textit{Remark.} Another way of representing $1$-cycle~$C$ as a sum of some cycles of length~$4$ is $\widehat \varphi \varphi C$ defined in our solution of~(\ref{en:knn:numcyc}).


(\ref{en:knn:numcyc}) \textit{Answer: $2^{(n-1)^2}$.} 

\textit{Remark.} Compare to Assertion~\ref{pr:orinumber}.\ref{en:orinumber:conngraph}.

\textit{Proof.}
It suffices to prove that there exists a 1--1 correspondence between $1$-cycles in $K_{n,n}$ and subsets\footnote{Sometimes the notation for a graph $K$ is used for the set of edges of $K$ (or, alternatively, the notation `$\cap K$' is used for the `restriction to $K$' operation on the sets of edges).}
of~$K_{n-1,n-1}$.
Define a map 
$$\varphi\colon H_1(K_{n,n}) \to 2^{K_{n-1,n-1}}\quad\text{by}\quad \varphi C := C \cap K_{n-1,n-1}.$$
For any edge $\sigma=ab'$ in $K_{n-1,n-1}$ denote by $\widehat \varphi\sigma$ the cycle $ab'nn'$ of length~$4$.
Define a map 
$$\widehat \varphi\colon 2^{K_{n-1,n-1}} \to H_1(K_{n,n})\quad\text{by}\quad 
\widehat \varphi D := \sum\limits_{\sigma \in D} \widehat \varphi\sigma.$$

Since $\varphi\widehat \varphi\sigma=\sigma$, we have $\varphi\widehat \varphi D = D$ for any $D \subset K_{n-1,n-1}$.
 
Take any $C \in H_1(K_{n,n})$. 
Since both $\widehat \varphi\varphi C$ and $C$ are $1$-cycles, $C' := \widehat\varphi\varphi C + C$ is a $1$-cycle. 
Applying $\varphi \widehat \varphi D = D$ for $D=\varphi C$, we obtain $\varphi C'=0$. 
So $C'\subset K_{n,n}\setminus K_{n-1,n-1}$. 
Hence by~\eqref{en:knn:tre} we have $C'=0$, so $\widehat \varphi\varphi C = C$.

Thus, $\varphi$ is a 1--1 correspondence.

(\ref{en:knn:rel}) See Assertion~\ref{pr:rook}.f.

{\it Remark.} 
Problem~\ref{p:knn} is equivalent to a particular case of Problem \ref{pr:rook} for $\ell=2$. 
An edge~$ab'\in K_{n,n}$ corresponds to $(a,b)\in [n]^2$.  
A vertex~$a\in[n]$ corresponds to the line $x=a$, 
a vertex~$b'\in[n]'$ corresponds to the line $y=b$.
An edge containing a vertex corresponds to an element of $[n]^2$ belonging to a line. 
Then $1$-cycles in~$K_{n,n}$ correspond to rook cycles in $[n]^2$. 
In particular, cycles of length~$4$ correspond to parallelepipeds (which for $\ell=2$ may be called parallelograms).


\smallskip
\textbf{\ref{pr:orinumber}.}
\emph{Answers:} (\ref{en:orinumber:tree}) $1$; \quad (\ref{en:orinumber:conngraph}) $2^{E-V+1}$. 

(\ref{en:orinumber:conngraph}) \emph{Proof.} \sophi{\phantomsection\label{en:sol:orinumber:conngraph}}
Let $T$ be a maximal tree of a connected graph $K$. 
It suffices to prove that there exists a 1--1 correspondence between $H_1(K)$ and subsets of $K \setminus T$.

Define a map 
$$\varphi\colon H_1(K) \to 2^{K \setminus T}\quad\text{by}\quad \varphi C := C \cap (K \setminus T).$$
For any edge $\sigma \in K \setminus T$ denote by $\widehat \varphi \sigma$ the simple cycle in $K$ formed by $\sigma$ and the simple path in $T$ joining ends of $\sigma$.
   
Define a map 
$$\widehat \varphi\colon 2^{K \setminus T} \to H_1(K)\quad\text{by}\quad 
\widehat \varphi D := \sum\limits_{\sigma \in D} \widehat \varphi \sigma.$$

Since $\varphi \widehat \varphi \sigma=\sigma$, we have $\varphi \widehat \varphi D = D$ for any $D\in 2^{K \setminus T}$.
 
Take any $C \in H_1(K)$. 
Since both $\widehat \varphi \varphi C$ and $C$ are $1$-cycles, $C':=\widehat \varphi \varphi C + C$ is a $1$-cycle. 
Applying $\varphi \widehat \varphi D = D$ for $D=\varphi C$, we obtain $\varphi C'=0$. 
So $C'\subset T$. 
Hence by~\eqref{en:orinumber:tree} we have $C'=0$, so $\widehat \varphi \varphi C = C$.
 
Thus, $\varphi$ is a 1--1 correspondence.

\smallskip
{\bf \ref{p:tkn}.} (\ref{en:tkn:4}) A~simple cycle in $\t{K_n}$ is formed by an $l$-element cyclic sequence of vertices such that 

$\bullet$  the number~$l$ is even and greater than~$2$; 

$\bullet$ the parts of vertices alternate in the sequence; 

$\bullet$ within each separate part, the vertices are not repeated; 

$\bullet$ there are no two consecutive vertices $m$ and $m'$.

Assume that such a cycle has no chords. 
Then any two non-consecutive vertices are not adjacent, i.e., they either are in the same part or are $m$ and $m'$. 
Any vertex has $\frac{l}{2}-2$ non-consecutive vertices from the other part. 
Then $\frac{l}{2}-2\le1$, so $l\in\{4,6\}$. 

\begin{figure}[h]\centering
\includegraphics[width=2cm]{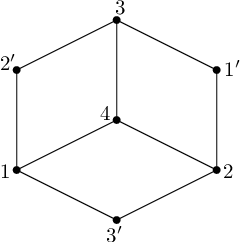}
\caption{A cycle of length 6 is the sum of three cycle of length 4}\label{f:cycles64}
\end{figure}

In a simple cycle of length~$6$ having no chords, any two opposite vertices are $m$ and $m'$. 
Hence the cycle is $(m_1,m'_2,m_3,m'_1,m_2,m'_3)$ for some pairwise distinct $m_1,m_2,m_3$. 
Since $n\ge4$, there is $a\in[n]-\{m_1,m_2,m_3\}$. 
Then the cycle equals 
$$(m_1,m_2',m_3,a')+(m_2,m_3',m_1,a')+(m_3,m_1',m_2,a'),$$ see Figure \ref{f:cycles64}. 
Now the required result holds by Assertion~\ref{pr:chl}.\ref{en:chl:nochord}.

(\ref{en:tkn:6},\ref{en:tkn:un}) Let $S$ be the set of all edges $ij'$ of $\t{K_n}$ such that $i,j>1$ and $(i,j)\ne(3,2)$.
The assertions follow (analogously to~Problem \ref{pr:orinumber}.\ref{en:orinumber:conngraph}) because 

$\bullet$ an edge $e'\in S$ is contained in the cycle $C_e$ if and only if $e=e'$; 

$\bullet$ the graph $\t{K_n}-S$ obtained from $\t{K_n}$ by deleting all edges of $S$ is a~tree. 
 
The second of these assertions follows because the edges of $\t{K_n}-S$ are exactly $32'$ and 
$1i',i1'$, where $i>1$, so $\t{K_n}-S$ is obtained by connecting with edge $32'$ the `star' trees consisting of edges $1i'$ with $i>1$, and of edges $1'i$ with $i>1$. 

\smallskip
{\bf \ref{p:symgra}.} (a,b) The natural projection $\t{K_n}\to K_n$ gives a 1--1 corespondence between $t$-symmetric 1-cycles in $\t{K_n}$ and 1-cycles in $K_n$. 
Use Assertions \ref{l:gencom1}.\ref{en:gencom1:nmb} and \ref{pr:chl}.\ref{en:chl:x-4}. 

(c) Take the base given by Assertions \ref{p:tkn}.\ref{en:tkn:6}\ref{en:tkn:un}.  
Then the 1-cycles 
$$\t{K_3}=C_{23'},\quad C_{ij'},\quad C_{23'}-C_{ji'},\quad\text{where}\quad i>j>1 \quad\text{and}\quad (i,j)\ne(3,2),$$  
also form a~base.
This base is as required because  $C_{23'}-C_{ji'}=1i'j1'32'12'31'23'=1i'j1'23'=tC_{ij'}$. 

Alternatively, one can take the cycles $\t{K_3}$, $1j'i2'$, $1'ji'2$, where $j>i>2$, and $1i'23'$, $1'i2'3$, where $i>3$.
 
\smallskip
{\bf \ref{3-homolzcx}.} (a) \textit{Statement.} There is only one integer 1-cycle in a tree. It assigns $0$ to all edges. 

\textit{Hint.} Start with proving that any integer $1$-cycle on a tree assigns~$0$ to all terminal edges.

(b) \textit{Statement.} Let $T$ be a maximal tree of a connected graph $K$. 
Every assignment of integers to edges of $K\setminus T$ has a unique extension to integer $1$-cycle in~$K$.

\textit{Proof.} Note that the sum of integer $1$-cycles is an integer $1$-cycle.

Denote by $\widehat \varphi_{\Z}(\sigma)$ the integer $1$-cycle in~$K$ 

$\bullet$ assigning $1$ to all edges of~$\widehat \varphi\sigma$ oriented in~$\widehat \varphi\sigma$ in the same way as~$\sigma$, 

$\bullet$ assigning~$-1$ to all edges of~$\widehat \varphi\sigma$ oriented in~$\widehat \varphi\sigma$ the opposite way, and 

$\bullet$ assigning~$0$ to all edges in~$K\setminus \widehat \varphi\sigma$.

For every assignment~$D$ of integers to edges of $K\setminus T$ denote $\widehat \varphi D:=\sum\limits_{\sigma\in K\setminus T} D (\sigma)\widehat \varphi_{\Z}(\sigma)$. 
Then $\widehat \varphi D$ is an integer $1$-cycle extending~$D$.

Take any integer $1$-cycle~$C$ in $K$ extending $D$. 
Then $C-\widehat \varphi D $ is an integer $1$-cycle that assigns~$0$'s to all edges in~$K\setminus T$. 
So by (a) integer $1$-cycle $C-\widehat \varphi D$ also assigns~$0$'s to all edges in~$T$. So $C=\widehat \varphi D$.


\smallskip
{\bf \ref{p:coh}.} (b) \textit{Answer:} $2^{E-V+1}$. 


\begin{figure}[!htb]
    \begin{minipage}[t]{0.33\textwidth}
	\centering
	\includegraphics[width=2cm]{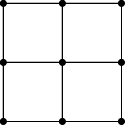}
    \end{minipage}
    \begin{minipage}[t]{0.33\textwidth}
	\centering
	\includegraphics[width=5cm]{real-33.mps}
    \end{minipage}
    \begin{minipage}[t]{0.33\textwidth}
	\centering
        \includegraphics[width=4cm]{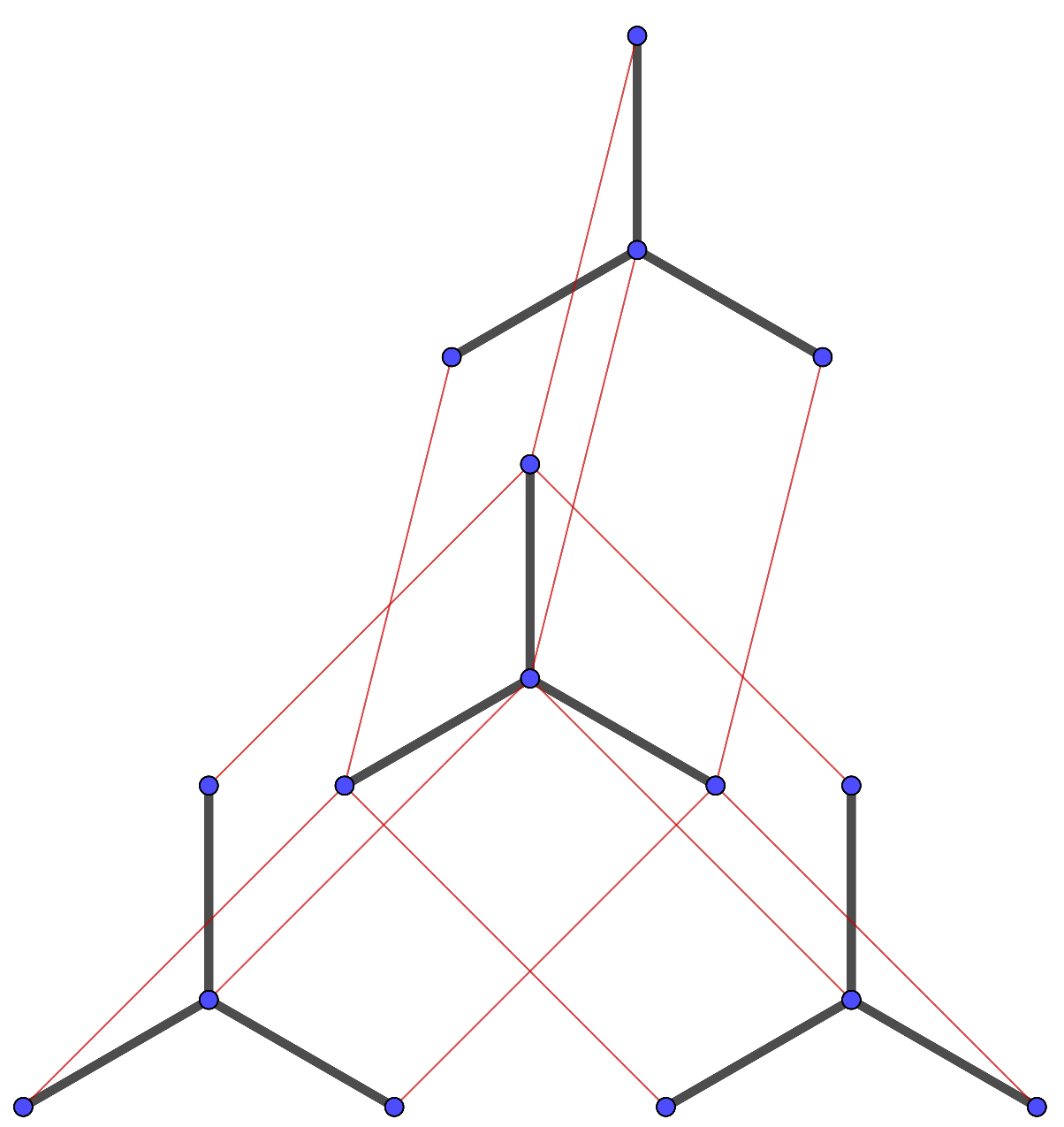}
    \end{minipage}\hfill
    \caption{Left: $K^{\square 2}_{2,1}$ is the $3\times 3$ grid graph. Middle: $K^{\square 2}_3$ on the torus $K^2_3$. Right: $K^{\square 2}_{3,1}$.}
        \label{f:grid-3-3}
        \label{f:k5i}
        \label{f:k3,1^2,2}
\end{figure}

\smallskip
{\bf\ref{hatg2}.}
(a) The $3\times 3$ grid graph, see Figure~\ref{f:grid-3-3}, left. 


	

(b) See Figure \ref{f:k5i}, middle (where boundary cycles are drawn in grey). 
Alternatively, $K_3^{\square 2}$ is obtained from the graph in (a) by adding edges between corresponding vertices:

$\bullet$ of the 1st and the 3rd rows; 

$\bullet$ of the 1st and the 3rd column.

(Compare to obtaining $K_3$ from the path~$123$ by adding an edge between $1$ and~$3$.)


(c) See Figure~\ref{f:k3,1^2,2}, right.

(d) Vertices are nodes of a~$4\times 4$ grid. 
Two vertices are joined by an edge if and only if they are in the same row or in the same column.
 
\smallskip
{\bf \ref{p:treek2}.} 
(b) The triodic cycle in $K_{3,1}^{\square 2}$ equals $\sum\limits_{i,j\in [3], i \not = j} i1' j1'$. 


(e,f) \textit{Answers:} no.  
 
(e) \textit{Hint.} The left projection of any boundary (and hence of any sum of boundaries) is empty.  

 
\smallskip
{\bf \ref{p:diag}.}
(a) We have $\diag(123)=1\times K_3+K_3\times 1+12\Box 23+12\Box 31 + 23\Box 31$.

(b) The off-diagonal cycle in $K_3^{\square 2}$ equals $1\times K_3+K_3\times 1+12\Box 31+31\Box 12$. 

(c) The antidiagonal cycle in $K_3^{\square 2}$ equals $1\times K_3+K_3\times1+23\Box31+31\Box23+31\Box31$. 

(d) We have  $1\times C+C\times 1 = \diag C+ \sum\limits_{i,j\in [n], i<j} i(i+1)\Box j(j+1)$, where $n+1:=1$.

(e) \textit{Answer:} no.  

\smallskip
{\bf \ref{p:k21}.}
\textit{Answers:} (a) $2^2$; (b) $2^2$; (c) $2^4$; (d) $2^6$.
 
\smallskip
{\bf \ref{stcycles2}.}
(a) Prove this by induction on the number of vertices of $T$, using deletion of a leaf vertex. 

(b) This can be proved by induction on the number of pairs of edges of $Z$ left-projected to the same edge of $K$. 
The inductive base when the number is zero can be proved by induction on the number of vertices~$b$ of~$K$ such that $Z\cap (b\times K)\ne\emptyset$.
For the inductive step, take such a pair $\{(\sigma,v_1), (\sigma,v_k)\}$, take a path $v_1\ldots v_k$ in $K$, and replace $Z$ by $Z+\sigma\Box v_1v_2+\ldots+\sigma\Box v_{k-1}v_k$.
 
(c) By (b), if $Z_x=Z_y=0$, then $Z\sim0$. 
Applying this to $Z=C + C_x\times a + a\times C_y$, we obtain the representation from (c). 
The uniqueness follows because the left and the right projections of sum of boundaries are equal to~$0$.

\begin{figure}[h]
\centerline{\includegraphics[width=4.5cm]{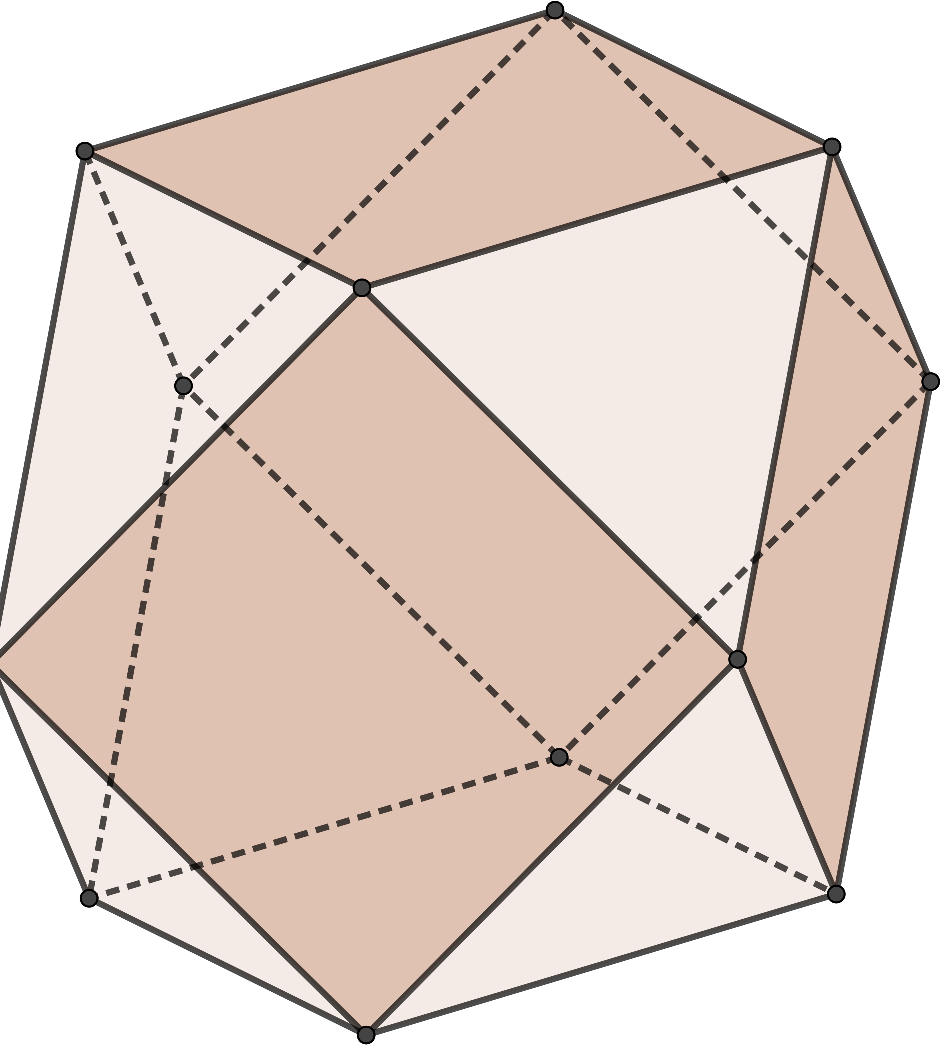}
\qquad\includegraphics[width=6cm]{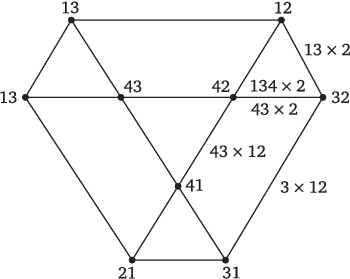}}
\caption{Left: cuboctahedron; the union of edges is $K_4^{\square\underline2}$ (three of the six boundaries are highlighted); the complement to triangular faces is $K_4^{\underline2}$. 
\newline
Right: the same with some explanations; the figure does not show the invisible part whose projection is obtained from the pictured projection by rotation through $\pi/3$
}
\label{f:del3}
\end{figure}
 
\smallskip
{\bf \ref{hatg}.} {\it Answers}.
(a) The disjoint union of two copies of $K_{2,1}$. 

(b) The cycle on 6 vertices.

(c) The cycle on 12 vertices.

(d) The union of edges of a cuboctahedron, see Figure~\ref{f:del3}.

\smallskip
{\bf \ref{p:knn21}.} {\it Answers}.

(a) $2^1$; (b) $2^5$; (c) $2^{17}$; (d) $2^{13}$; (e) $2^{43}$; (f) $2^{41}$.

(a') $2^1$; (b') $2^1$; (c') $2^{5}$; (d') $2^{7}$; (e') $2^{8}$; (f') $2^{12}$.

{\it Hint.} Apply \ref{pr:orinumber}.b to obtain (a)-(f).
 
{\it Remark.} Cf. examples before Problem \ref{t:dpcycle}. 

\smallskip
{\bf \ref{p:stcycles}.} (a) 
It suffices to prove the result for $K$ connected. 
It suffices to prove that if a sum of pairwise distinct boundaries is a symmetric 1-cycle, then either the sum is zero, or there are two summands symmetric to each other. 
The analogous result for $K^{\square 2}$ replaced by $K\Box T\bigcup\limits_{T\Box T} T\Box K$, where $T$ is a tree, is proved by induction on the number of vertices in $T$, using deletion of a leaf vertex. 
For reduction of (a) to this analogous result, use 
the 
equality $\sum\limits_{i,j} \sigma_i\Box\tau_j = 0$ for two simple cycles in $K$ having  consecutive edges $\sigma_1\ldots\sigma_k$ and $\tau_1\ldots\tau_\ell$.   
Using this equality we replace the given sum of boundaries by an equal sum not containing summands $\sigma\Box\tau$ corresponding to edges $\sigma,\tau$ outside a maximal tree of $K$. 


\smallskip
\textbf{\ref{l:gencom}.} 
(\ref{en:gencom:gen})
It suffices to prove that an arbitraty $2$-cycle $C$ equals the sum $\widehat C := \sum\limits_{\{ i, j, k\} \in C,\ i,j,k<n} T_{i,j,k,n}$ of tetrahedra $T_{i,j,k,n}$ over $3$-element subsets $\{ i, j, k\} \in C$ not containing $n$. 
This follows because every face in $ C + \widehat C $ contains $n$, so by~(\ref{en:gencom:tre}) we have $C + \widehat C = 0$, i.e., $C=\widehat C$.

 
(\ref{en:gencom:rel})
A rigorous formulation is as follows (we identify tetrahedra in $[n]$ with 4-element subsets of $[n]$). 
Suppose that $R$ is a set of $4$-element subsets of $[n]$ such that any $3$-element subset of $[n]$ is contained in an even number of those $4$-element subsets that are elements of $R$. 
Then there is a family $R'$ of $5$-element subsets of $[n]$ such that a $4$-element subset $B$ of $[n]$ is contained in $R$ if and only if $B$ is contained in an odd number of those $5$-element subsets that are elements of $R'$.  

This statement is a higher dimensional analogue of Assertion \ref{l:gencom}.\ref{en:gencom:gen}; the proof is analogous.
 
\smallskip
\textbf{\ref{pr:rook}.} 
(c) For $a\in[n-1]^\ell$ denote by $P(a):=\{n,a_1\}\times\ldots\times\{n,a_\ell\}$ the parallelepiped with opposite vertices $a$ and $(n,\ldots,n)$. 
It suffices to prove that any rook cycle $C\subset[n]^\ell$ equals the sum $\widehat C$ of parallelepipeds $P(a)$ over $a\in C\cap[n-1]^\ell$. 
The sum $C+\widehat C$ is a rook cycle. 
Since $P(a)\cap[n-1]^\ell=\{a\}$, we have $(C+\widehat C)\cap[n-1]^\ell=\emptyset$.
By (b), we have $C+\widehat C = \emptyset$, i.e., $C=\widehat C$.

(f) \emph{Hint.} Any parallelepiped~$P$ with some $P_i\subset[n-1]$ is a sum of some given relations and two parallelepipeds obtained from~$P$ by replacing with~$n$ one of the two elements in~$P_i$.  
So in each of the two parallelepipeds the number of pairs $P_j$ not containing $n$ is smaller than in~$P$. 
Hence 
any relation between parallelepipeds is a sum of some given relations and a~relation $P(a_1)+\ldots+P(a_s)=0$ for some pairwise distinct $a_1,\ldots,a_s\in[n-1]^{\ell}$.
The latter relation is trivial because 
\[
    \emptyset=(P(a_1)+\ldots+P(a_s))\cap[n-1]^{\ell}=\{a_1,\ldots,a_s\}.
\]
Here the second equality holds because $P(a)\cap[n-1]^{\ell}=\{a\}$.


\smallskip
\textbf{\ref{p:realcyl}.} 
Let $O,V,A_{11},\ldots A_{1n}$ be points in $\R^3$ of which no 4  lie in the same  plane.
For every $p\in[n]$ take a point $A_{2p}$ such that $\vv{OA_{2p}}=\vv{OV}+\vv{OA_{1p}}$.
If $V$ is close enough to $O$, then the points $A_{jp}$, $j\in\{1,2\}$, $p\in[n]$, are as required, 
i.e., the union of parallelograms $A_{1p}A_{2p}A_{2q}A_{1q}$ forms the drawing of the cylinder over $K_n$ in $\R^3$.


    \begin{figure}[h]
    \centerline{\includegraphics[width=5cm]{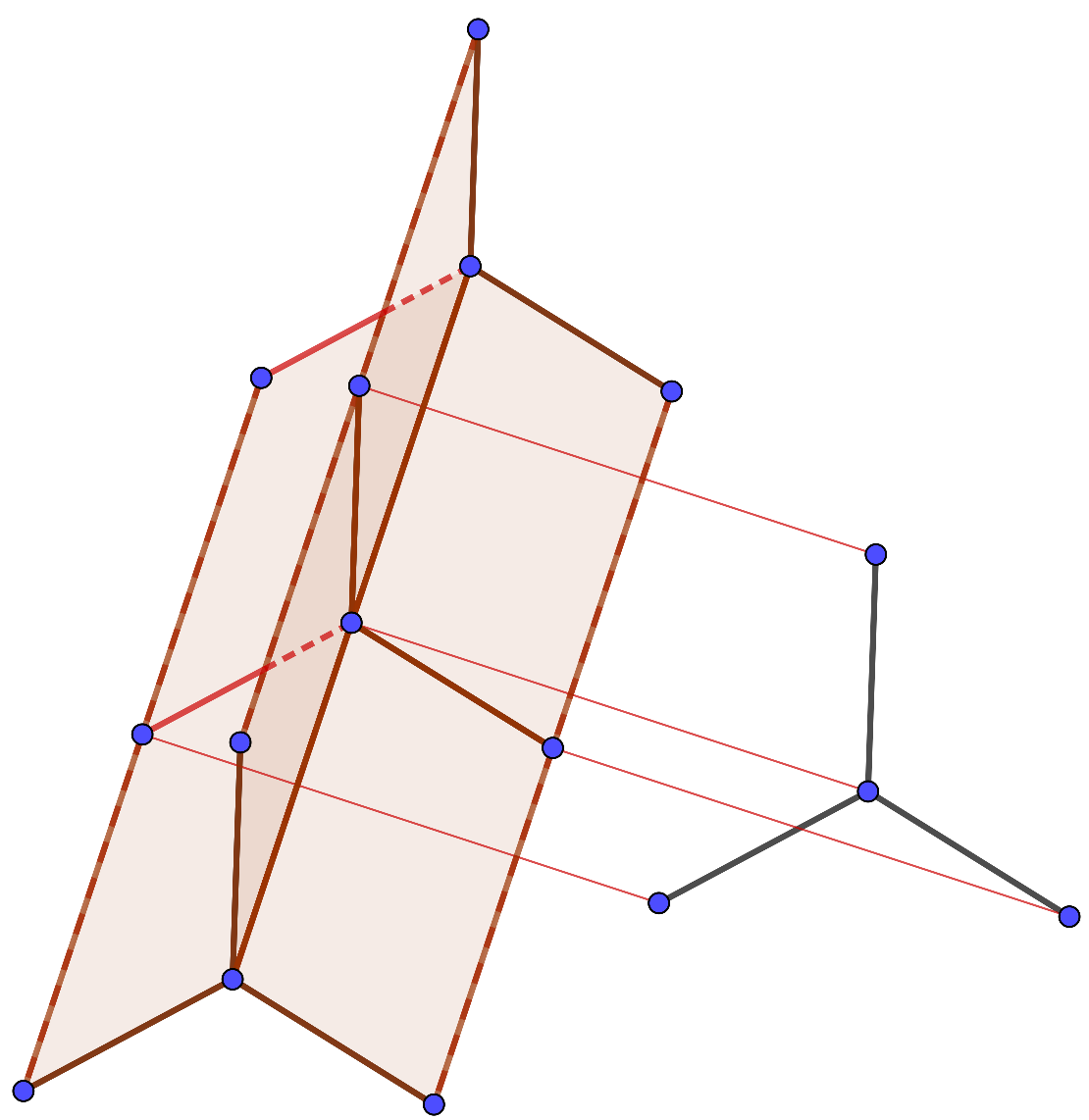}}
    \caption{To drawing of $K_{3,1}^2$ in 3-space}
    \label{f:k3,1^2,1}
    \end{figure}


\smallskip
\textbf{\ref{p:cyl}.} 
(b) See Figure \ref{f:k3,1^2,1}. 

\smallskip
\textbf{\ref{p:ce2cy}.}
(b) In this proof by $\sigma, \tau$ we denote edges of $K$, and by $v$ we denote a vertex of $K$: 
$$C \text{ is a 2-cycle} \ \Leftrightarrow 
\ \left\{ 
\begin{array}{cc}
     & |\{ \sigma  :\ \sigma\ni v,\ (\sigma, \tau) \in C \}| \text{ is even for any $v$ and } \tau \\
     & |\{ \tau  :\ \tau\ni v,\ (\sigma, \tau) \in C \}| \text{ is even for any $v$ and } \sigma
\end{array} 
\right. \ \Leftrightarrow$$
$$\Leftrightarrow 
\ \left\{ 
\begin{array}{cc}
     & |\{ \sigma \in C_{\cdot, \tau}\ :\ \sigma\ni v\}| \text{ is even for any $v$ and } \tau\\
     & |\{ \tau \in C_{\sigma, \cdot}\ :\ \sigma\ni v\}| \text{ is even for any $v$ and } \sigma
\end{array} 
\right.
\ \Leftrightarrow 
\ \left\{ 
\begin{array}{cc}
     & C_{\cdot, \tau} \text{ is a 1-cycle for any } \tau\\
     & C_{\sigma, \cdot} \text{ is a 1-cycle for any } \sigma.
\end{array} 
\right.
$$

\smallskip
\textbf{\ref{t:torus}.}
(a) Let $C := \widetilde C_1 \times \widetilde C_2 \subset K^2$ be the product of simple cycles $\widetilde C_1$ and $\widetilde C_2$ in $K$. 
For every edge $\sigma$ of $K$ the set $C_{\sigma,\cdot}$ is either $\widetilde C_2$ or the empty set. 
For every edge $\tau$ of $K$ the set $C_{\cdot, \tau}$ is either $\widetilde C_1$ or the empty set. 
Hence by~Assertion \ref{p:ce2cy}.b, $C$ is a cellular $2$-cycle.


(c) \emph{Answer:} no.

\emph{Solution.}
Let $C$ be a cellular $2$-cycle in $T^2$. 
Since $T$ is a tree, by Assertion \ref{p:ce2cy}.b $C_{\sigma,\cdot}=0$ for any edge $\sigma$ of $T$. 
Hence $C=0$. 

(d) \emph{Answer:} no.

\emph{Solution.}
Let $C$ be a cellular $2$-cycle in $\overline T$.
Since $T$ is a tree, by~Assertion \ref{p:ce2cy}.b $C_{\sigma,\cdot}=0$ for any edge $\sigma$ of $K \setminus T$. 
Hence $C\subset T\times K$. 
Since $T$ is a tree, by~Assertion \ref{p:ce2cy}.b $C_{\cdot,\sigma}=0$ for any edge $\sigma$ of $K$. 
Hence $C=0$. 

\smallskip
\textbf{\ref{p:k2}.}
\emph{Answers:}
(a) 2; \quad (b) 2; \quad (c) $2^4$; \quad (d) $2^9$. 

\emph{Solutions.}

(a, b) Let $C$ be a cellular $2$-cycle in $K^2$ such that $C \neq K$, i.e., $(\sigma, \tau) \notin C$ for some edges $\sigma, \tau$ in $K$. 
By Assertion \ref{p:ce2cy}.b, $C_{\sigma, \cdot}$ is a $1$-cycle. 
Since $K$ is a cycle graph, $C_{\sigma, \cdot} = 0$. 
Then $(\sigma, \sigma) \notin C$. 
Hence $C \subset \overline{K - \sigma}$. 
By answer to Problem~\ref{t:torus}.d, $C = 0$.

(c) By answer to Problem \ref{t:kunn}.d, the number of $2$-cycles in $K_{2,3}^2$ equals $2^{(E - V + 1)^2} = 2^{2^2} = 16$.

(d) By answer to Problem \ref{t:kunn}.d, the number of $2$-cycles in $K_{4}^2$ equals $2^{(E - V + 1)^2} = 2^{3^2} = 2^9 = 512$.

\smallskip
\textbf{\ref{t:dpcycle}.}
(a) \emph{Answer:} no.

\emph{Solution.}
Let $C$ be a cellular $2$-cycle in $K^{\underline 2}$. 
Take an edge $\sigma$ of $K$. 
By Assertion \ref{p:ce2cy}.b, $C_{\sigma, \cdot}$ is a $1$-cycle. 
Since $C_{\sigma, \cdot}$ is a subset of the path graph $K - \sigma$, then $C_{\sigma, \cdot} = 0$. 
Hence $C = 0$.

(b) \emph{Answer:} no.

\emph{Solution.}
Let $C$ be a cellular $2$-cycle in $K^{\underline 2}$. 
By Assertion \ref{p:ce2cy}.b, $C_{\sigma, \cdot},\ C_{\cdot, \sigma}$ are $1$-cycles for any edge $\sigma$ of $K$. 
Since  both $C_{0j, \cdot}$ and $C_{\cdot, 0j}$ are subsets of the path graph $K - 0 - j$ for any $j \in [n]$, then $C_{0j, \cdot} = C_{\cdot, 0j} = 0$ for any $j \in [n]$. 
Then $C\subset(K - 0)^{\underline 2}$ for the cycle graph $K - 0$. 
Hence by (a), $C = 0$. 

(c) Denote $C := K^{\underline 2}_{3, 3}$. 
For any edge $\sigma$ in $K_{3, 3}$ both $C_{\sigma, \cdot}$ and $C_{\cdot, \sigma}$ are cycles of length~4.
Then the required result follows by Assertion \ref{p:ce2cy}.b.

(d) Denote $C := K^{\underline 2}_5$. 
For any edge $\sigma$ in $K_5$ both $C_{\sigma, \cdot}$ and $C_{\cdot, \sigma}$ are cycles of length~3.
Then the required result follows by Assertion \ref{p:ce2cy}.b.

(e) Take $K := K_6$ and $C := 123 \times 456$. 
Since $(12, 34) \notin C$, the product $C$ is a non-empty proper subset of $K^{\underline 2}$. 
By Assertion \ref{t:torus}.a, the product $C$ of cycles is a cellular $2$-cycle. 
 
\smallskip
\textbf{\ref{t:kunn}.}
(a, b) \emph{Hint:} the required representations can be obtained from the proof of (c).

(c) Part (c) easily follows from the particular case of connected $K$, which we further assume.

Let $T$ be a maximal tree of $K$.
Define a map 
$$\varphi\colon H_2(K^2) \to 2^{K^2-\overline{T}} \quad\text{by}\quad \varphi C := C \cap (K^2 - \overline{T}).$$
For any edge $\sigma \in K \setminus T$ the $1$-cycle $\widehat \varphi \sigma$ is defined in the proof of Assertion~\ref{pr:orinumber}.\ref{en:orinumber:conngraph}.
Define a map $$\widehat \varphi\colon 2^{K^2-\overline{T}} \to H_2(K^2) \quad\text{by}\quad \widehat \varphi D := \sum\limits_{(\sigma, \tau) \in D} \widehat \varphi \sigma \times \widehat \varphi \tau.$$

By Assertion~\ref{t:torus}.d, $\widehat \varphi \varphi C = C$ for any $C \in H_2(K^2)$, i.e., an arbitraty cellular $2$-cycle $C$ in $K^2$ is the sum of tori $\widehat \varphi \sigma \times \widehat \varphi \tau$ over $(\sigma, \tau) \in \varphi C$.

(d) \emph{Answer:} $2^{(E-V+1)^2}.$

\emph{Solution.} In (c), we defined maps $\varphi: H_2(K^2) \to 2^{K^2-\overline{T}}$ and $\widehat \varphi: 2^{K^2-\overline{T}} \to H_2(K^2)$ for a maximal tree $T$ of connected graph $K$ and proved that $\widehat \varphi \varphi C = C$ for any $C \in H_2(K^2)$. 
In addition, $\varphi \widehat \varphi D = D$ for any $D \in 2^{K^2-\overline{T}}$. 
Thus, $\varphi$ is a 1--1 correspondence.  
Hence 
$$|H_2(K^2)| = |2^{K^2-\overline{T}}| = |2^{(K \setminus T)^2}| = 2^{(E-V+1)^2}.$$

\smallskip
\textbf{\ref{t:kunneth}.}
(b) 
Denote $B := \{C_i \times C_j\ :\ i, j \in [q] \}$.
In the following two paragraphs we prove that every $2$-cycle in $K^2$ is a sum of some $2$-cycles from $B$, and that such a representation is unique, respectively.

Apply Assertion \ref{t:kunn}.c. 
Every torus $\widetilde C_1 \times \widetilde C_2 \subset K^2$ is in turn a sum of some $2$-cycles from $B$, because both $\widetilde C_1$ and $\widetilde C_2$ are sums of some $1$-cycles from $C_1, \ldots, C_q$. 

Denote by $N$ the number of connected components of $K$. 
Analogously to Assertion \ref{pr:orinumber}.\ref{en:orinumber:conngraph}, $|H_1(K)| = 2^{E - V + N}$. 
Since $C_1,\ldots, C_q$ is a base in $H_1(K)$, it follows that $2^q = |H_1(K)|$, so $q = E - V + N$. 
Therefore $|B| = q^2 = (E - V + N)^2$. 
By generalization of the answer to Problem \ref{t:kunn}.d to an arbitrary graph $K$, we have $|H_2(K^2)| = 2^{(E - V + N)^2}$, which   
is the number of linear combinations of $2$-cycles from $B$. 
This implies the uniqueness.

\smallskip
\textbf{\ref{t:kunnethsym}.}
(a) \emph{Hint:} this follows by Assertion \ref{p:ce2cy}.b.

(b) \emph{Answer:} no.

\emph{Solution.}
Take $K := K_5$ (or $K := K_{3,3}$) and $C := K^{\underline2}$. 
By Assertion \ref{t:dpcycle}.c (or \ref{t:dpcycle}.d), $C$ is a cellular $2$-cycle. 
Since there are no two non-empty vertex-disjoint cycles in $K$, the 2-cycle $C$ is not a sum of some products of vertex-disjoint cycles.

\smallskip
\textbf{\ref{t:2cyc-bij}.} 
(\ref{en:2cyc-bij:num}) \textit{Answer:} $2^{(n^2-3n+1)^2}$.


(\ref{en:2cyc-bij:compl})
\emph{Hint.} The 1--1 correspondence is defined by the formula $f(\sigma_1\sigma_2',\tau_1\tau_2') := (\sigma_1\tau_1',\sigma_2\tau_2')$ for disjoint edges $\sigma,\tau\in[n]^2$.
(Recall that for $a,b\in[n]$ we denote by $ab'$ the edge joining $a$ and $b'$.) 
Cells $\alpha, \beta \in [n] \times [n]' \times [n] \times [n]'$ in $K_{n,n}^{2}$ are adjacent if and only if they differ only at one of their four coordinates. Since $K_{n,n}^{\underline2}, \widetilde{K_{n}}^{2} \subset K_{n,n}^{2}$, it follows that $\alpha$ and $\beta$ in $K_{n,n}^{\underline2}$ are adjacent if and only if $f\alpha$ and $f\beta$ in $\widetilde{K_{n}}^{2}$ are adjacent, i.e., $f$ respects the adjacence.

{\it Remark.}
Define a 1--1 correspondence $\widehat f$ between edges of $K_{n,n}^{\square \underline2}$ and edges of $\t{K_n}^{\square 2}$ by skipping one of $\sigma_1,\sigma_2',\tau_1,\tau_2'$ in the formula for $f$ (e.g. $\widehat f(\sigma_1,\tau_1\tau_2') = (\sigma_1\tau_1',\tau_2')$).
One can check that the correspondences $f,\widehat f$ respect incidence.
(Recall that an edge $\eta$ and a $2$-cell $\alpha=(\xi,\zeta)$ are \emph{incident} if $\eta=\xi$ or $\eta=\zeta$.)
Note that there is no 1--1 correspondence between vertices that saves incidence.
Indeed, there is no 1--1 correspondence between $1$-cycles in $K_{n,n}^{\square \underline2}$ and in $\t{K_n}^{\square 2}$.
E.g. for $n=3$ we have $|H_1(\t{K_3}^{\square 2})| = 2^{37}$, as opposed to the answer to Problem~\ref{p:knn21}.e;
besides, in $\t{K_3}^{\square 2}$ there are four $1$-cycles up to boundaries, as opposed to the answer to Problem~\ref{p:knn21}.e'.


\smallskip
\textbf{\ref{p:tsym}.}
(c) \emph{Answer:} no.

\emph{Solution.} Take $K=K_3$ and $C=K^2$. 
Then 
$C$ is a non-empty symmetric cellular $2$-cycle in $K^2$. 
Since $K\times K + K \times K = 0$, there are no non-empty symmetrized tori in $K^2$. 
Hence $C$ is not a sum of some symmetrized tori.

\smallskip
\textbf{\ref{l:h2sym}.}
(a) For any symmetric cellular $2$-cycle $C$ in $K^{\underline2}$ we have
$$C \stackrel{(1)}=  \sum_{(\sigma,\tau)\in C-\overline T}\widehat \varphi \sigma\times\widehat \varphi \tau
        \stackrel{(2)}= \sum_{\{\sigma,\tau\}\colon(\sigma,\tau)\in C-\overline T}
        (\widehat \varphi \sigma\times\widehat \varphi \tau + \widehat \varphi \tau\times\widehat \varphi \sigma). $$
Here 

$\bullet$ (1) holds by the formula in the hint to Assertion \ref{t:kunn}.c, and 

$\bullet$ (2) holds because $C-\overline T$ is symmetric, and $C-\overline T \subset C \subset K^{\underline2}$, so 
$\sigma \neq \tau$ for any pair $(\sigma, \tau) \in C-\overline T$.






(c,d,e) See the details in \cite{SS23, DMN+}.


{\it Books, surveys and expository papers in this list are marked by the stars.}

\end{document}